\documentclass[11pt]{amsart}

\usepackage{amsmath,amssymb,latexsym,soul,cite,mathrsfs}
\usepackage{tikz}
\usepackage{color,enumitem,graphicx}
\usepackage[colorlinks=true,urlcolor=blue,
citecolor=red,linkcolor=blue,linktocpage,pdfpagelabels,
bookmarksnumbered,bookmarksopen]{hyperref}
\usepackage[english]{babel}

\usepackage[left=3.3cm,right=3.3cm,top=2.8cm,bottom=2.8cm]{geometry}
\usepackage[hyperpageref]{backref}

\usepackage[colorinlistoftodos]{todonotes}
\makeatletter
\providecommand\@dotsep{5}
\def\listtodoname{List of Todos}
\def\listoftodos{\@starttoc{tdo}\listtodoname}
\makeatother

\usepackage[notref,notcite]{showkeys}

\numberwithin{equation}{section}

\newtheorem{theorem}{Theorem}[section]

\newtheorem{proposition}[theorem]{Proposition}

\begin{document}

\title[Fractional $p$-Laplacian equations with sandwich pairs]{Fractional $p$-Laplacian equations with sandwich pairs}

\author{J. Vanterler da C. Sousa}

\address[J. Vanterler da C. Sousa]
{\newline\indent Aerospace Engineering, PPGEA-UEMA
\newline\indent
Department of Mathematics, DEMATI-UEMA
\newline\indent
São Luís, MA 65054, Brazil.}
\email{\href{vanterler@ime.unicamp.br}{vanterler@ime.unicamp.br}}

\pretolerance10000


\begin{abstract} The main purpose of this paper is to consider new sandwich pairs and investigate the existence of solution for a new class of fractional differential equations with $p$-Laplacian via variational methods in $\psi$-fractional space $\mathbb{H}^{\alpha,\beta;\psi}_{p}(\Omega)$. The results obtained in this paper are the first to make use of the theory of $\psi$-Hilfer fractional operators with $p$-Laplacian.

\end{abstract}

\maketitle
\section{Introduction and motivation}

In this paper we consider a new class of fractional differential equations with $p$-Laplacian given by
\begin{equation}\label{1.3}
_{C}^{\bf H}\mathfrak{D}_{T}^{\alpha ,\beta ;\psi }\left( \left\vert ^{\bf H}\mathfrak{D}_{0+}^{\alpha
,\beta ;\psi }\phi\right\vert ^{p-2}\text{ }^{\bf H}\mathfrak{D}_{0+}^{\alpha ,\beta ;\psi
}\phi\right) =f(\xi,\phi),\text{ in }\Omega
\end{equation}
where $\Omega =[0,T]$ is a bounded domain in $\mathbb{R}$ (where $\mathbb{R}$ is the real line), $_{C}^{\bf H}\mathfrak{D}_{T}^{\alpha ,\beta ;\psi }(\cdot )$ and $^{\bf H}\mathfrak{D}_{0+}^{\alpha ,\beta ;\psi }(\cdot )$ are Hilfer-Caputo and $\psi$-Hilfer fractional derivative
of order $\frac{1}{p}<\alpha <1$ and type $0\leq \beta \leq 1$, $1<p<\infty $ and $f$ is a Caratheodory function on $[0,T]\times \mathbb{R} $ with subcritical growth. 
Note that, for $p>2$ the (\ref{1.3}) degenerate and $1< p < 2$ singular, at points where $\left\vert^{\bf H}\mathfrak{D}_{0+}^{\alpha
,\beta ;\psi }\phi\right\vert=0$. For $p=2$ we just get the usual, that is, $_{C}^{\bf H}\mathfrak{D}_{T}^{\alpha ,\beta ;\psi }\left(\text{ }^{\bf H}\mathfrak{D}_{0+}^{\alpha ,\beta ;\psi}\phi\right)$. I this sense, we put the condition $\left\vert^{\bf H}\mathfrak{D}_{0+}^{\alpha
,\beta ;\psi }\phi\right\vert\neq 0$.

The theory of fractional differential equations, over the last few years, has attracted attention, from problems involving a theoretical approach, to problems involving applications \cite{Samko,Kilbas}. We highlight problems of controllability theory of differential equation solutions \cite{Vijayakumar,Vijayakumar1,Vijayakumar2,Vijayakumar3,Vijayakumar4}, which are relevant problems that contribute significantly to the area.

First, before commenting on sandwich pairs, it is worth mentioning the importance and relevance of the p-Laplacian equations and their contributions. Over the years, the impact and importance of p-Laplacian equations to the theory of differential equations, especially the elliptical ones, is remarkable and undeniable. Problems of existence and multiplicity are indeed interesting and have drawn attention over these decades, especially in this last decade, with double phase problems \cite{novo1,novo2,novo3,novo4}. We can highlight some applications, such as: in the mechanics of nano structures, fluids, diffusion process and asymptotic dynamics \cite{1,2,11,14,18,19,20} and the references therein. On the other hand, it is worth highlighting the problems of differential equations with $p$-Laplacian via variational methods and $\psi$-Hilfer fractional operators, which began in mid 2021 and has been gaining ground in the area \cite{Sousa1, Sousa80,Sousa,Sousa33,parte1,parte2,parte3,parte4,parte5}. The theory is still new and recent, and under construction. In this sense, results in the area are still quite restricted. Consequently, there are two aspects, i.e., the first is that there are few results to use and, in the vast majority of cases, it is necessary to build them. On the other hand, this allows a range of options to work with and numerous open problems arise as the theory grows.

The first ideas of sandwich pairs were built using eigenspaces and used to find critical points of a functional. This approach was taken by Schechter \cite{6}. We can also highlight the works on solved quasilinear problems using cones as sandwich pairs \cite{4,5}. See to the works \cite{7,10,chen}.

Many problems arising in science and engineering call for the solving of the Euler equations of functionals, i.e., equations
of the form $G'(u)=0$ where $G(u)$ is a $C^{1}$ functional arising from the given data. Since the development of the calculus of variations there has been interest in finding critical points of functionals. This was intensified by the fact that for many equations arising in practice, the solutions are critical points. See some interesting papers on sandwich pairs and applications \cite{par1,par2,par3,par4}.

Perera and Scheter \cite{4}, discussed boundary value problem for the $p$-Laplacian using the notion of sandwich pairs, that is, they addressed the following problem
\begin{equation}\label{op23}
\left\{ 
\begin{array}{ccc}
-\Delta_{p} u & = & f(x,u),\text{ in }\Omega \\ 
u & = & 0,\,\, on\,\partial \Omega.
\end{array}
\right. 
\end{equation}
For more details about the problem (\ref{op23}), see \cite{4}.

The sandwich pairs used until 2007, were introduced using the eigenspaces of a semilinear operator and are therefore unsuitable for dealing with quasilinear problems where there are no eigenspaces. In this sense, in 2008 Perera and Scheter \cite{5}, showed that the method could be modified to be applied in the problem (\ref{op23}).

Perera and Schechter \cite{Schechter} discussed the solution of problems of type
\begin{equation}\label{op}
\left\{ 
\begin{array}{ccc}
-\Delta_{p} u & = & \nabla_{u}  F(x,u),\text{ in }\Omega \\ 
u & = & 0,\,\, on\,\partial \Omega.
\end{array}
\right. 
\end{equation}
For more details about the problem (\ref{op}), see \cite{Schechter}.

We say shall that a pair of subsets $\mathcal{A},\mathcal{B}$ of a Banach space $W$ forms a sandwich pair, if for any $\mathcal{E}^{\alpha ,\beta ;\psi }(\cdot )\in
C^{1}(W,\mathbb{R})$ the inequality \cite{par4}
\begin{equation}\label{1.1}
-\infty <b:=\underset{\mathcal{B}}{\inf }\,\,\mathcal{E}^{\alpha ,\beta ;\psi }(\cdot )\leq a:= \underset{\mathcal{A}}{\sup }\,\,\mathcal{E}^{\alpha ,\beta ;\psi }(\cdot )<\infty
\end{equation}
implies that there is a sequence $(\phi_{j})\subset W$ such that 
\begin{equation}\label{1.2}
\mathcal{E}^{\alpha ,\beta ;\psi }(\phi_{j})\rightarrow c,\left( \mathcal{E}^{\alpha ,\beta ;\psi
}\right) ^{\prime }(\phi_{j})\rightarrow 0
\end{equation}
for some $c\in \lbrack b,a]$.

Note that the sequence satisfying (\ref{1.1}) is called a Palais–Smale sequence at the level $c$ and $\mathcal{E}^{\alpha ,\beta ;\psi }$ satisfies the compactness condition $(PS)_{c}$ if every such sequence
has a convergent subsequence.

Motivated by the works \cite{4,5,Schechter}, in this paper we concern to attack the existence of a solution to the problem (\ref{1.3}) via variational methods and sandwich pairs. We will discuss the existence of solutions through two theorems, one with the lower limit condition and the other with the upper limit condition and from the function
\begin{equation*}
\Theta(\xi,t)=\mathfrak{F}(\xi,t)-t f(\xi,t),
\end{equation*}
where $\mathfrak{F}(\xi,t)=\displaystyle\int_{0}^{t}f(x,s )ds$, in other words, we are interested in discussing the following results:

\begin{theorem}\label{Theorem1.2} If 
\begin{equation}\label{1.8}
(\lambda _{l}+\varepsilon )|t|^{p}-\mathcal{V}(\xi)\leq  \mathfrak{F}(\xi,t)\leq \lambda
_{l+1}|t|^{p}+\mathcal{V}(\xi)
\end{equation}
for some $l,\varepsilon >0$ and $\mathcal{V}\in L_{\psi }^{1}(\Omega )$, and 
\begin{equation}\label{1.9}    
\Theta(\xi,t)\leq C(|t|+1),\,\overline{\Theta}(\xi):=\underset{|t|\rightarrow \infty }{
\overline{\lim }}\frac{\Theta(\xi,t)}{|t|}<0,\,\,a.e
\end{equation}
then the problem {\rm(\ref{1.3})} has a solution.
\end{theorem}

\begin{theorem}\label{Theorem1.3} If 
\begin{equation}\label{1.10}
\lambda _{l}|t|^{p}-\mathcal{V}(\xi)\leq  \mathfrak{F}(\xi,t)\leq (\lambda _{l+1}-\varepsilon
)|t|^{p}+\mathcal{V}(\xi)
\end{equation}
for some $l,\varepsilon >0$ and $\mathcal{V}\in L_{\psi }^{1}(\Omega )$, and 
\begin{equation}\label{1.11}
    \Theta(\xi,t)\geq -C(|t|+1),\,\,\underline{\Theta}(\xi):=\underset{|t|\rightarrow \infty }{%
\underline{\lim }}\frac{\Theta(\xi,t)}{|t|}>0,\,\,a.e
\end{equation}
then the problem {\rm(\ref{1.3})} has a solution.
\end{theorem}

A natural consequence when working with fractional operators is to obtain the classic case, as a particular case, this is of paramount importance and relevance for the investigated results. Here in this work, it is possible to obtain such a property, in addition to obtaining other possible particular cases from the choice of the parameters $\beta \rightarrow 1$ or $\beta \rightarrow 0$ and, from the function $\psi(\cdot)$. The discussion of some cases will be discussed at the end of the paper, as "special cases". However, one of the limitations of this work is that it cannot choose the function $\psi(\xi)=\ln\,\xi$ as a particular case, since the problem (\ref{1.3}) is being covered in the space $\psi$-fractional $\mathbb{H}^{\alpha,\beta;\psi}_{p}(\Omega)$. However, it is possible to discuss this case, but it is necessary to work with the weight space of $\mathbb{H}^{\alpha,\beta;\psi}_{p}(\Omega)$. Furthermore, we can rule out that the results obtained here are the first in the area of fractional differential equations with $p$-Laplacian and $\psi$-Hilfer fractional operators. Certainly, the results presented in this work will draw attention to future work, in particular, a natural continuation of this work as highlighted at the end of the paper.

In section 2, we present definitions and results about the theory of fractional operators and sandwich pairs. Finally, in Section 3, we investigate the main results of the article, i.e., the proof of {\bf Theorem \ref{Theorem1.2}} and {\bf Theorem \ref{Theorem1.3}}. In this sense, we present some special special cases.
\section{Mathematical background: preliminaries}

In this section, we present some definitions of $p$-integrable spaces, $\psi$-fractional space, and results on the $\psi$-Hilfer fractional derivative. Also for, sandwich pair results, we end the section.

The space of $p$-integrable function with respect to a function $\psi$ is defined as
\begin{eqnarray}\label{ppp}
    L^{p}_{\psi}([a,b],\mathbb{R})=\left\{\phi:[a,b]\rightarrow\mathbb{R}: \int_{a}^{b} |\phi(x)|^{p} \psi'(x) dx<\infty \right\}
\end{eqnarray}
with norm
\begin{equation*}
    ||\phi||_{L^{p}_{\psi}([a,b],\mathbb{R})} = \left(\int_{a}^{b} \psi'(x) |\phi(x)|^{p} dx \right)^{1/p}.
\end{equation*}

Choosing $p=1$ in Eq.(\ref{ppp}), we have the integrable space $L^{1}_{\psi}([a,b],\mathbb{R})$ with its respective norm 
\begin{equation*}
    ||\phi||_{L^{1}_{\psi}([a,b],\mathbb{R})} = \int_{a}^{b} \psi'(x) |\phi(x)| dx.
\end{equation*}

Let $n-1<\alpha<n$ with $n\in\mathbb{N}$, $I=[a,b]$ is the interval such that $-\infty\leq a<b\leq \infty$ and there exist two functions $f,\psi \in C^{n}([a,b],\mathbb{R})$ such that $\psi$ increasing and $\psi'(\xi)\neq 0$, for all $\xi\in I$. The $\psi$-Hilfer fractional derivatives left-sided and right-sided ${^{\mathbf H}\mathfrak{D}}^{\alpha,\beta;\psi}_{a}(\cdot)$ $({^{\mathbf H}\mathfrak{D}}^{\alpha,\beta;\psi}_{b}(\cdot))$ of order $\alpha$ and type $0\leq \beta \leq 1$ are defined by \cite{J1}
\begin{equation}\label{derivada1}
{^{\mathbf H}\mathfrak{D}}^{\alpha,\beta;\psi}_{a}\phi(\xi)= {\bf I}^{\beta(n-\alpha),\psi}_{a} \left(\frac{1}{\psi'(\xi)} \frac{d}{d\xi} \right) {\bf I}^{(1-\beta)(n-\alpha),\psi}_{a} \phi(\xi)
\end{equation}
and
\begin{equation}\label{derivada2}
{^{\mathbf H}\mathfrak{D}}^{\alpha,\beta;\psi}_{b}\phi(\xi)= {\bf I}^{\beta(n-\alpha),\psi}_{b} \left(-\frac{1}{\psi'(\xi)} \frac{d}{d\xi} \right) {\bf I}^{(1-\beta)(n-\alpha),\psi}_{b} \phi(\xi)
\end{equation}
where
\begin{equation*}
    {\bf I}^{\alpha,\psi}_{a} \phi(\xi)=\dfrac{1}{\Gamma(\alpha)} \int_{a}^{\xi} \psi'(s)(\psi(\xi)- \psi(s))^{\alpha-1} \phi(s) ds,\,\,to\,\,a<s<\xi
\end{equation*}
and
\begin{equation*}
    {\bf I}^{\alpha,\psi}_{b} \phi(\xi)=\dfrac{1}{\Gamma(\alpha)} \int_{\xi}^{b} \psi'(s)(\psi(s)- \psi(\xi))^{\alpha-1} \phi(s) ds,\,\,to\,\,\xi<s<b
\end{equation*}
are the fractional integrals of $\phi$ with respect to the function $\psi$. The definition of $_{C}^{\bf H}\mathfrak{D}_{T}^{\alpha ,\beta ;\psi } (\cdot)$ is the commutation between the integral operators $ {\bf I}^{\beta(n-\alpha),\psi}_{a}(\cdot)$ and ${\bf I}^{(1-\beta)(n-\alpha),\psi}_{a}(\cdot) $ of Eq.(\ref{derivada2}). Furthermore, $\left(\dfrac{d}{d\xi} \right)$ and $\psi'(\xi)$ are classical derivatives. 

Note that the function $\psi(\cdot)$ is part of kernel of the $\psi$-Riemann-Liouville fractional integral and, consequently, of the $\psi$-Hilfer fractional derivative. The motivation to introduce the $\psi$-Hilfer fractional derivative comes from the fractional derivatives: Caputo, Riemann-Liouville, and Hilfer, to unify a wide class of fractional operators in a unique operator. The restriction on the function $\psi(\xi)$, is that $\psi(\xi)\neq 0$ for all $\xi \in I$, since $\dfrac{1}{\psi' (\xi)}$, as stated in the definition itself. A priori, a physical meaning about $\psi(\cdot)$, is still unknown. However, what can be discussed are their respective particular cases, based on the choice of the function $\psi(\cdot)$ and the limits of $\beta\rightarrow 1$ or $\rightarrow 0$. For example, in the particular choice of $\psi(\xi)=\xi$, $\psi(\xi)=\xi^{\rho}$ with $\rho>0$ and $\psi(\xi) =\ln (\xi)$, classic fractional derivatives are obtained, for example Caputo, Riemann-Liouville, Hadamard, Caputo-Hadamard, Hilfer, among others.

Let $\alpha\in (0,1)$, $p\in (1,\infty)$ and $\frac{1}{p}+\frac{1}{q}\leq 1+\alpha$. If $\phi\in L^{q}_{\psi}([a,b])$ and $\varphi\in L^{p}_{\psi}([a,b])$, then the following integration by parts \cite{Sousa1}
\begin{equation}\label{eq.218}
\int_{a}^{b}\left( {\bf I}_{a}^{\alpha ;\psi }\varphi\left( \xi\right) \right) \phi\left( \xi\right) \psi'(\xi) d\xi =\int_{a}^{b} \varphi\left( \xi\right) {\bf I}_{b}^{\alpha ;\psi }\phi\left( \xi\right) \psi ^{\prime }\left( \xi\right) d\xi.
\end{equation}

On the other hand, let $\alpha\in (0,1)$, $\beta\in [0,1]$, $\phi,\varphi\in AC[a,b]$ and $-\dfrac{1}{\psi'(\cdot)} {\bf I}_{b}^{\beta(1-\alpha);\psi}\phi\in L^{2}_{\psi}[a,b]$, then \cite{Sousa}
\begin{eqnarray}\label{A}
&&\int_{a}^{b}\,^{{\bf H}}_{c}\mathfrak{D}_{b}^{\alpha,\beta ;\psi } \phi(\xi) \varphi(\xi)\psi'(\xi) d\xi\notag\\&&=-\lim_{x\rightarrow b} {\bf I}^{(1-\alpha)(1-\beta);\psi}_{0+}\varphi(\xi) \,\,{\bf I}^{\beta(1-\alpha);\psi}_{b}\phi(\xi)\notag\\
&&+ \lim_{x\rightarrow 0+} {\bf I}^{(1-\alpha)(1-\beta);\psi}_{0+}\varphi(\xi) {\bf I}^{\beta(1-\alpha);\psi}_{b}\phi(\xi) + \int_{0}^{b} \phi(\xi) ^{{\bf H}}\mathfrak{D}_{0+}^{\alpha,\beta ;\psi }  \varphi(\xi) \psi'(\xi) d\xi.
\end{eqnarray}

The $\psi$-fractional space is given by \cite{Sousa1}
\begin{equation*}
\mathbb{H}^{\alpha,\beta;\psi}_{p}(\Omega)=\left\{\phi\in L^{p}_{\psi}(\Omega)~:~\left\vert^{\rm H}\mathfrak{D}^{\alpha,\beta;\,\psi}_{0+}\phi\right\vert\in L^{p}_{\psi}(\Omega),\,\, \phi=0\,\,on\,\, \partial\Omega\right\}
\end{equation*}
with the norm 
\begin{equation*}
||\phi||=||\phi||_{\mathbb{H}^{\alpha,\beta;\psi}_{p}(\Omega)}=||\phi||_{L^{p}_{\psi}(\Omega)}+||^{\rm H}\mathfrak{D}^{\alpha,\beta;\psi}_{0+}\phi||_{L^{p}_{\psi}(\Omega)}.
\end{equation*}
The space $C_{0}^{\infty }(\Omega )$ is dense in $\mathbb{H}_{p}^{\alpha ,\beta ;\psi }(\Omega )$. The space $\mathbb{H}^{\alpha,\beta;\psi}_{p}(\Omega)$ is separable and reflexive Banach space \cite{Sousa1}.

Indeed, for $\varphi \in C_{0}^{\infty }([0,T])$ and taking the integral in both sides of Eq.(\ref{1.3}), yields
\begin{equation}\label{*}
\int_{0}^{T}\text{ }_{C}^{\bf H}\mathfrak{D}_{T}^{\alpha ,\beta ;\psi }\left( \left\vert
^{\bf H}\mathfrak{D}_{0+}^{\alpha ,\beta ;\psi }\phi\right\vert ^{p-2}\text{ }%
^{\bf H}\mathfrak{D}_{0+}^{\alpha ,\beta ;\psi }\phi\right) \varphi (\xi)\psi ^{\prime
}(\xi)d\xi=\int_{0}^{T}f(\xi,\phi)\varphi (\xi)\psi ^{\prime }(\xi)d\xi.
\end{equation}

Using the relation Eq.(\ref{A}) and taking 
\begin{eqnarray*}
    \underset{%
x\rightarrow 0+}{\lim }{\bf I}^{(1-\alpha )(1-\beta );\psi }\varphi (\xi)=0=\underset {x\rightarrow T}{\lim }{\bf I}^{(1-\alpha )(1-\beta );\psi }\varphi (\xi),
\end{eqnarray*}
we have
\begin{eqnarray}\label{**}
&&\int_{0}^{T}\text{ }_{C}^{\bf H}\mathfrak{D}_{T}^{\alpha ,\beta ;\psi }\left( \left\vert
^{\bf H}\mathfrak{D}_{0+}^{\alpha ,\beta ;\psi }\phi\right\vert ^{p-2}\text{ }%
^{\bf H}\mathfrak{D}_{0+}^{\alpha ,\beta ;\psi }\phi\right) \varphi (\xi)\psi ^{\prime
}(\xi)d\xi\notag\\&&=\int_{0}^{T}\left\vert ^{\bf H}\mathfrak{D}_{0+}^{\alpha ,\beta ;\psi }\phi\right\vert
^{p-2}\text{ }^{\bf H}\mathfrak{D}_{0+}^{\alpha ,\beta ;\psi }\phi\text{ }^{\bf H}\mathfrak{D}_{0+}^{\alpha
,\beta ;\psi }\varphi (\xi)\psi ^{\prime }(\xi)d\xi.
\end{eqnarray}

Therefore, from Eq.(\ref{*}) and Eq.(\ref{**}), yields
\begin{equation}
\int_{0}^{T}\left\vert ^{\bf H}\mathfrak{D}_{0+}^{\alpha ,\beta ;\psi }\phi\right\vert ^{p-2}%
\text{ }^{\bf H}\mathfrak{D}_{0+}^{\alpha ,\beta ;\psi }\phi\text{ }^{\bf H}\mathfrak{D}_{0+}^{\alpha ,\beta
;\psi }\varphi (\xi)\psi ^{\prime }(\xi)d\xi=\int_{0}^{T}f(\xi,\phi)\varphi (\xi)\psi
^{\prime }(\xi)d\xi.
\end{equation}

Consider $\varphi =\phi$, we have
\begin{equation}
\int_{0}^{T}\left\vert ^{\bf H}\mathfrak{D}_{0+}^{\alpha ,\beta ;\psi }\phi\right\vert ^{p}
\text{ }\psi ^{\prime }(\xi)d\xi=\int_{0}^{T}f(\xi,\phi)\phi\psi ^{\prime }(\xi)d\xi.
\end{equation}

Now, we define the Euler functional $\mathcal{E}^{\alpha ,\beta ;\psi }:\mathbb{H}_{p}^{\alpha
,\beta ;\psi }([0,T],\mathbb{R})\rightarrow \mathbb{R}$ on $\mathbb{H}_{p}^{\alpha ,\beta ;\psi }([0,T],\mathbb{R})$, given by
\begin{equation}\label{1.4}
\mathcal{E}^{\alpha ,\beta ;\psi }(\phi):=\frac{1}{p}\int_{0}^{T}\left\vert
^{\bf H}\mathfrak{D}_{0+}^{\alpha ,\beta ;\psi }\phi\right\vert ^{p}\text{ }\psi ^{\prime
}(\xi)d\xi-\int_{0}^{T}\mathfrak{F}(\xi,\phi(\xi))\psi ^{\prime }(\xi)d\xi
\end{equation}
where $\mathfrak{F}(\xi,t)=\displaystyle\int_{0}^{t}f(x,s )ds$.

The solution of (\ref{1.3}) coincide with the critical points of the $C^{1}$ functional Eq.(\ref{1.4}).

Consider the result about sandwich pairs.
\begin{proposition}{\rm\cite{4,Schechter}}\label{Proposition2.1} $\Gamma $ be the class of maps $\gamma \in C\left( W\times
\lbrack 0,1],W\right) $ such that:

$(a) \gamma _{0}=id$;

$(b) \underset{(\phi,t)\in W\times \lbrack 0,1]}{\sup }\left\Vert \gamma
_{t}(\phi)-\phi\right\Vert <\infty $
where $\gamma _{t}=\gamma (\cdot ,t)$. Assume that for any $\gamma \in
\Gamma $,
\begin{equation}\label{2.1}
\gamma (\mathcal{A})\cap \mathcal{B}\neq \emptyset.
\end{equation}
Then $\mathcal{A},\mathcal{B}$ forms a sandwich pair.    
\end{proposition}

Let $\mathcal{A}_{1}\neq \emptyset$ and $\mathcal{B}_{1}\neq \emptyset $ be subsets of the interval $J$ in a Banach space $W$ such that $dist(\mathcal{A}_{1},\mathcal{B}_{1})>0$. We say that $\mathcal{A}_{1}$ links $ \mathcal{B}_{1}$ if for any $\mathcal{E}^{\alpha ,\beta ;\psi }(\cdot )\in C^{1}(J,\mathbb{R})$
\begin{equation}\label{2.2}
-\infty <\underset{\mathcal{A}_{1}}{\sup }\,\,\mathcal{E}^{\alpha ,\beta ;\psi }(\cdot
)=:a_{0}<b_{0}:=\underset{\mathcal{B}_{1}}{\inf }\,\,\mathcal{E}^{\alpha ,\beta ;\psi }(\cdot
)<\infty
\end{equation}
implies that there exists a sequence $(\phi_{j})\subset J$ such that Eq.(\ref{2.2}) holds for some $c\geq b_{0}$.

\begin{proposition}{\rm\cite{9}}\label{Proposition2.2} $\mathcal{A}_{1}$ links $\mathcal{B}_{1}$ in $J$ if for any $\varphi \in
C(C\mathcal{A}_{1},J)$ such that $\varphi (\cdot ,0)=id_{\mathcal{A}_{1}}$ 
\begin{equation}\label{2.3}
\varphi (C\mathcal{A}_{1})\cap \mathcal{B}_{1}\neq \emptyset 
\end{equation}
where $C\mathcal{A}_{1}=(\mathcal{A}_{1}\times \lbrack 0,1])/(\mathcal{A}_{1}\times \{1\})$ is a subset on $\mathcal{A}_{1}$.
\end{proposition}

\begin{proposition}{\rm\cite{4,Schechter}}\label{Proposition2.3} If $\mathcal{A}_{1}$ and $\mathcal{B}_{1}$ satisfy the hypotheses of {\bf Proposition \ref{Proposition2.2}} in $J$, then
\begin{equation}\label{2.4}
\mathcal{A}=\pi ^{-1}(\mathcal{A}_{1})\cup \{0\},\text{ }\mathcal{B}=\pi ^{-1}(\mathcal{B}_{1})\cup \{0\}
\end{equation}
forms a sandwich pair, where $\pi :W\backslash \{0\}\rightarrow J$.
\end{proposition}

Consider the nonlinear eigenvalue fractional problem
\begin{equation}\label{1.5}
\left\{ 
\begin{array}{ccc}
_{C}^{\bf H}\mathfrak{D}_{T}^{\alpha ,\beta ;\psi }\left( \left\vert ^{\bf H}\mathfrak{D}_{0+}^{\alpha
,\beta ;\psi }\phi\right\vert ^{p-2}\text{ }^{\bf H}\mathfrak{D}_{0+}^{\alpha ,\beta ;\psi
}\phi\right) & = & \lambda |\phi|^{p-2},\text{ in }\Omega \\ 
\phi & = & 0,\,\,on\,\,\partial\Omega.
\end{array}%
\right. 
\end{equation}

It's eigenvalues are similar the critical values of the $C^{1}$ functional
\begin{equation}\label{1.6}
\mathcal{I}_{\psi}(\phi)=\frac{1}{\displaystyle\int_{0}^{T}\psi ^{\prime }(\xi)|\phi|^{p}d\xi}
\end{equation}
on the interval $J$ in $\mathbb{H}_{p}^{\alpha ,\beta ;\psi }([0,T],\mathbb{R})$. Let $\Gamma ^{l}$ be the class of odd continuous maps $\gamma $ from the interval $J_{l-1}$ in $\mathbb{R}$ to $J$ and set
\begin{equation}\label{1.7}
\lambda _{l}=\underset{\gamma \in \Gamma _{l}}{\inf }\,\,{\max _{\phi\in \gamma
(J_{l-1})}}\,\,\mathcal{I}_{\psi}(\phi).
\end{equation}

Then $0<\lambda _{1}<\lambda _{2}\leq \cdot \cdot \cdot \rightarrow \infty $
are eigenvalues of the problem (\ref{1.5}). Consider
\begin{equation*}
\Theta(\xi,t)=\mathfrak{F}(\xi,t)-t f(\xi,t).
\end{equation*}

We emphasize that the resonance is considered only concerning the specific variational eigenvalues given by (\ref{1.7}) and not with respect to other possible non-variational eigenvalues or variational eigenvalues which are given by different methods.
\section{Main results}

This section, we concern to investigate the main results of this paper, i.e., the proof of {\bf Theorem \ref{Theorem1.2}} and {\bf Theorem \ref{Theorem1.3}} through the results presented in Section 2.

So, we start with the proof of the first result according to the theorem below:

\begin{proof} {\bf (Proof of Theorem \ref{Theorem1.2})} Using Eq.(\ref{1.7}) there is a $\gamma\in \Gamma_{l}$ such that $\mathcal{I}_{\psi}\leq \lambda_{l}+\dfrac{\varepsilon}{2}$ on $\mathcal{A}_{1}=\gamma(J_{l-1})$. Let $\mathcal{B}_{1}=\left\{ \phi\in J:\mathcal{I}_{\psi}(\phi)\right\}\geq \lambda_{l+1}$. Since $\lambda_{l}+\dfrac{\varepsilon}{2}<\lambda_{l+1}$ by inequality (\ref{1.8}), $\mathcal{A}_{1}$ and $\mathcal{B}_{1}$ are disjoint. Since $\mathcal{A}_{1}$ is compact and $\mathcal{B}_{1}$ is closed, it follows that $dist(\mathcal{A}_{1},\mathcal{B}_{1})>0$. We claim that $\mathcal{A}_{1}$ links $\mathcal{B}_{1}$ in $J$. Given $\varphi\in C(C\mathcal{A}_{1}, J)$ such that $\varphi(\cdot,0)=id_{\mathcal{A}_{1}}$, writing $\xi\in J_{l}$ as $(\xi',\xi_{l+1})\in \mathbb{R}\oplus\mathbb{R}$, define $\overline{\gamma}\in\Gamma_{l+1}$ by
\begin{equation*}
\overline{\gamma}(\xi)=\left\{ 
\begin{array}{ccc}
\varphi(\gamma(\xi'/|x'|),\xi_{l+1}) \,\,for\,\,0\leq \xi_{l+1}<1 \\ 
\varphi(\mathcal{A}_{1}\times \left\{1\right\}) \,\,for\,\,\xi_{l+1}=1  \\ 
\overline{\gamma}(\xi',-\xi_{l+1}) \,\,for\,\,\xi_{l+1}<0  \\ 
\end{array}.
\right. 
\end{equation*}

Then, $\overline{\gamma} (J_{l})\cap \mathcal{B}_{1}\neq \emptyset$ by definition of $\lambda_{l+1}$ so Eq.(\ref{2.3}) holds as $\mathcal{B}_{1}$ is symmetric. Hence $\mathcal{A},\mathcal{B}$ given by (\ref{2.4}) forms a sandwich pair by {\bf Proposition \ref{Proposition2.3}}. Let $\mathcal{E}^{\alpha,\beta;\psi}$ given by Eq.(\ref{1.4}). Since
\begin{equation*}
    \int_{0}^{T}\psi'(\xi) \left\vert^{\bf H}\mathfrak{D}_{0+}^{\alpha ,\beta ;\psi }\phi\right\vert ^{p} d\xi\geq \lambda_{l+1} \int_{0}^{T} \psi'(\xi) |\phi|^{p} d\xi,\,\,\phi\in \mathcal{B}
\end{equation*}
and \begin{equation*}
    \int_{0}^{T} \psi'(\xi) \left\vert ^{\bf H}\mathfrak{D}_{0+}^{\alpha ,\beta ;\psi }\phi\right\vert ^{p} d\xi\leq (\lambda_{l+1}+\varepsilon) \int_{0}^{T} \psi'(\xi) |\phi|^{p} d\xi,\,\,\phi\in \mathcal{A},
\end{equation*}
(\ref{1.8}) implies
\begin{equation*}
    -\int_{0}^{T} \psi'(0) \mathcal{V}(\xi) d\xi\leq \inf_{\mathcal{B}} \mathcal{E}^{\alpha,\beta;\psi}(\cdot)\leq \sup_{\mathcal{A}} \mathcal{E}^{\alpha,\beta;\psi}\leq \int_{0}^{T} \psi'(\xi) \mathcal{V}(\xi) d\xi.
\end{equation*}
Hence there exists a sequence $(\phi_{j})\subset \mathbb{H}^{\alpha,\beta;\psi}_{p}(\Omega)$ satisfying Eq.(\ref{1.2}). 

Since $(\phi_{j})$ is bounded and so there exists a convergent subsequence $\rho_{j}=||\phi_{j}||\rightarrow\infty$, a subsequence of $\overline{\phi}_{j}=\dfrac{\phi_{j}}{\rho_{j}}\rightharpoonup \overline{\phi}$ in $\mathbb{H}^{\alpha,\beta;\psi}_{p}(\Omega)$, strongly in $L_{\psi}^{p}(\Omega)$ and almost everywhere in $\Omega$. Then, using Eq.(\ref{1.2}), yields
\begin{eqnarray*}
    \int_{0}^{T} \psi'(\xi) \frac{\Theta(\xi,\phi_{j})}{\rho_{j}} d\xi &&= \int_{0}^{T} \psi'(\xi) \left(\frac{\mathfrak{F}(\xi,\phi_{j})- tf(\xi,\phi_{j})}{\rho_{j}} \right) d\xi\notag\\
    &&=\int_{0}^{T}\psi'(\xi) \frac{\mathfrak{F}(\xi,\phi_{j})}{\rho_{j}} d\xi  \int_{0}^{T}\psi'(\xi) \frac{t f(\xi,\phi_{j})}{\rho_{j}} d\xi \notag\\
    &&=\dfrac{\dfrac{ (\mathcal{E}^{\alpha,\beta;\psi}(\phi_{j}))' \phi_{j}}{p^{2}} - \mathcal{E}^{\alpha,\beta;\psi}(\phi_{j})}{\rho_{j}} \rightarrow 0.
\end{eqnarray*}

Using $\overline{\phi}_{j}=\dfrac{\phi_{j}}{\rho_{j}}$ and inequality (\ref{1.9}), it's follows that
\begin{equation*}
    \overline{\lim} \int_{0}^{T} \psi'(\xi) \frac{\Theta(\xi,\phi_{j})}{\rho_{j}}d\xi\notag\\
    \leq \int_{0}^{T} \psi'(\xi) \overline{\lim} \frac{\Theta(\xi,\phi_{j})}{|\phi_{j}|}  |\overline{\phi}_{j}|d\xi\notag\\
    = \int_{0}^{T} \psi'(\xi) \overline{\Theta}(\xi) |\overline{\phi}|d\xi \leq 0.
\end{equation*}

Since $\overline{\Theta}<0$ almost everywhere, we have that $\overline{\phi}=0$ almost everywhere. In this sense, passing to the limit, yields
\begin{equation*}
    1-\frac{\mathcal{E}^{\alpha,\beta;\psi}(\phi_{j})}{\rho_{j}^{p}} = \int_{0}^{T} \psi'(\xi) \frac{\mathfrak{F}(\xi,\phi_{j})}{\rho_{j}^{p}} d\xi \leq \int_{0}^{T} \psi'(\xi) \left(\lambda_{l+1} |\overline{\phi_{j}}|^{p}+\frac{\mathcal{V}}{\rho_{j}^{p}} \right) d\xi
\end{equation*}
gives a contraction.    
\end{proof}

Now, we will prove the second main result of this paper.

\begin{proof} {\bf (Proof of Theorem \ref{Theorem1.3})} Let a sequence $(\varepsilon_{j})\subset (0,\varepsilon]$ decreasing to $0$ and 
\begin{equation*}
    \mathcal{E}^{\alpha,\beta;\psi}_{j} (\phi)= \mathcal{E}^{\alpha,\beta;\psi}(\phi) - \varepsilon_{j} \int_{0}^{T} \psi'(\xi) |\phi|^{p} d\xi.
\end{equation*}
Then, using Eq.(\ref{1.10}), we obtain
\begin{equation}
    \left((\lambda_{l}+\varepsilon_{j})|t|^{p}-\mathcal{V}(\xi)\right)\leq \left(\mathfrak{F}(\xi,t)+\varepsilon_{j}|t|^{p}\right) \leq \left(\lambda_{l+1}|t|^{p}+\mathcal{V}(\xi)\right).
\end{equation}

In this sense, there exists a sequence $(\phi_{j})\subset \mathbb{H}^{\alpha,\beta;\psi}_{p}(\Omega)$ such that $ \mathcal{E}^{\alpha,\beta;\psi}_{j}(\phi_{j})$ is bounded $ (\mathcal{E}^{\alpha,\beta;\psi}_{j})'(\phi_{j})\rightarrow 0$ (see proof of {\bf Theorem \ref{Theorem1.2}}). Since that $(\phi_{j})$ is bounded and so there exists a subsequence converges to a critical point of $ \mathcal{E}^{\alpha,\beta;\psi}_{j}(\cdot)$. If $\rho_{j}=||\phi_{j}||_{}\psi\rightarrow\infty$, a subsequence of $\overline{\phi}_{j}=\dfrac{\phi_{j}}{\rho_{j}}\rightharpoonup\overline{\phi}$ in $\mathbb{H}^{\alpha,\beta;\psi}_{p}(\Omega)$, strongly in $L^{p}_{\psi}(\Omega)$, and almost everywhere $\Omega=[0,T]$. Then using Eq.(\ref{1.10}), yields
\begin{eqnarray*}
    \int_{0}^{T} \psi'(\xi) \frac{\Theta(\xi,\phi_{j})}{\rho_{j}} d\xi &&= \int_{0}^{T} \left( \frac{\mathfrak{F}(\xi,\phi_{j})-tf(\xi,\phi_{j})}{\rho_{j}}\right)d\xi\notag\\
    &&= \int_{0}^{T} \psi'(\xi) \frac{\mathfrak{F}(\xi,\phi_{j})}{\rho_{j}}d\xi -\int_{0}^{T} \psi'(\xi) \frac{tf(\xi,\phi_{j})}{\rho_{j}}d\xi \notag\\
    &&=\dfrac{\dfrac{ (\mathcal{E}^{\alpha,\beta;\psi}(\phi_{j}))' \phi_{j}}{p^{2}} - \mathcal{E}^{\alpha,\beta;\psi}(\phi_{j})}{\rho_{j}} \rightarrow 0.
\end{eqnarray*}

On the other hand, using Eq.(\ref{1.11}) we get
\begin{eqnarray*}
    &&\underline{\lim}\int_{0}^{T} \psi'(\xi) \frac{\Theta(\xi,\phi_{j},\phi_{j})}{\rho_{j}} d\xi\notag\\
   && \geq \int_{0}^{T} \psi'(\xi) \underline{\lim} \frac{\Theta(\xi,\phi_{j})}{|\phi_{j}|} |\overline{\phi}_{j}| d\xi
    \notag\\
    &&\geq  \int_{0}^{T}\psi'(\xi) \underline{\Theta}(\xi) |\overline{\phi}|\geq 0.
\end{eqnarray*}

Since $\underline{\Theta}>0$ almost everywhere then $\overline{\phi}=0$ almost everywhere. In this sense,  passing to the limit in
\begin{equation*}
    1-\frac{\mathcal{E}^{\alpha,\beta;\psi}(\phi_{j})}{\rho_{j}^{p}} = \int_{0}^{T} \psi'(\xi) \left(\frac{\mathfrak{F}(\xi,\phi_{j})}{\rho_{j}^{p}}+\varepsilon_{j} |\phi_{j}|^{p} \right) d\xi \leq \int_{0}^{T} \psi'(\xi) \left(\lambda_{l+1} |\overline{\phi_{j}}|^{p}+\frac{\mathcal{V}}{\rho_{j}^{p}} \right) d\xi
\end{equation*}
given a contradiction. Therefore, we concluded the proof.
\end{proof}

\subsection{Special cases}

A natural consequence of the results obtained above is the freedom to present a wide class of possible particular cases, especially the integer case. In this sense, we present three cases below:

{\bf Case 1:} (Special case.) Taking the limit $\alpha\rightarrow 1$ and $\psi(\xi)=\xi$, we get the problem in its classic version, given by
\begin{equation}\label{klo}
\left( \left\vert \phi'\right\vert ^{p-2}\,\,\phi'\right)' =f(\xi,\phi),\text{ in }\Omega.
\end{equation}

{\bf Case 2:} (Caputo fractional operator case.) Taking $\beta=1$ and $\psi(\xi)=\xi$, we have
\begin{equation}\label{klo1}
_{C}\mathfrak{D}_{T}^{\alpha}\left( \left\vert ^{\bf C}\mathfrak{D}_{0+}^{\alpha}\phi\right\vert ^{p-2}\text{ }^{\bf C}\mathfrak{D}_{0+}^{\alpha
}\phi\right) =f(\xi,\phi),\text{ in }\Omega.
\end{equation}

{\bf Case 3:} (Riemann-Liouville fractional operator case.) Taking $\beta=0$ and $\psi(\xi)=\xi$, we get
\begin{equation}\label{klo2}
_{C}^{\bf RL}\mathfrak{D}_{T}^{\alpha}\left( \left\vert ^{\bf RL}\mathfrak{D}_{0+}^{\alpha}\phi\right\vert ^{p-2}\text{ }^{\bf RL}\mathfrak{D}_{0+}^{\alpha}\phi\right) =f(\xi,\phi),\text{ in }\Omega.
\end{equation}

Once we obtain the particular cases, i.e., the problems (\ref{klo})-(\ref{klo2}), the results investigated here, {\bf Theorem \ref{Theorem1.2}} and {\bf Theorem \ref{Theorem1.3}} , are valid for such case. It is possible to notice that the freedom of choice of the function $\psi(\cdot)$, allows to obtain other examples involving fractional operators, however, we restrict ourselves to the ones previously discussed. See other formulations of fractional operators that can be obtained by choosing the function $\psi(\cdot)$ \cite{J1}.
\section{Conclusion and future work}

At the end of this work, we were able to obtain the existence of a solution for a new class of fractional differential equations with $p$-Laplacian and sandwich pairs via variational methods. Although the results investigated here are new, there is still a long way to go to be investigated, which elucidates that the theory is still under construction and, consequently, some open problems a priori, we can highlight as follows:
\begin{enumerate}
    \item Taking the problem (\ref{1.3}), we can work with the $p(x)$-Laplacian has a more complex nonlinearity that raises some of the essential difficulties, for example, it is inhomogeneous.

    \item An interesting issue that can also be worked on, is to discuss the problem (\ref{1.3}) with double phase.

    \item Finally, a possible discussion of the problem (\ref{1.3}) with Kirchhoff problem would be interesting.

\end{enumerate}

In this sense, we conclude this work in the certainty that the work paves the way for new results.

\vspace{1.5em}

\section*{Acknowledgements}

J. Vanterler da C. Sousa very grateful to the anonymous reviewers for their useful comments that led to improvement of the manuscript.

{\bf Funding} Funding information is not applicable/No funding was received.
\vspace{1.5em}

{\bf Availability of data and materials} Data sharing not applicable to this paper as no data sets were generated or analyzed during the current study.

\section*{Declarations}

{\bf Conflict of interest} The author have no conflicts to disclose.
\vspace{1.5em}

{\bf Ethical approval} Not applicable.



\end{document}